\documentclass[12pt]{amsart}
\usepackage{amssymb}

\def\cc{{\mathcal C}}
\def\dist{D^*}
\def\emp{{\emptyset}}
\def\ll{{\mathcal L}}
\def\nn{{\mathbb N}}
\def\pp{{\mathcal P}}
\def\qq{{\mathbb Q}}
\def\reals{{\mathbb R}}
\def\rmand{\mbox{ and }}
\def\rmiff{\mbox{ iff }}
\def\sq{\subseteq}
\def\st{{\;:\;}}
\def\text#1{\makebox(0,0)[cc]{#1}}
\def\th{{\theta}}
\def\zz{{\mathbb Z}}

\newtheorem{theorem}{Theorem}[section]
\newtheorem{corollary}[theorem]{Corollary} 

\newtheorem{prop}[theorem]{Proposition}
\newtheorem{question}[theorem]{Question}
\theoremstyle{definition}
\newtheorem{definition}[theorem]{Definition}
\newtheorem{example}[theorem]{Example}

\begin{document}

\title{Steinhaus Sets and Jackson Sets}
\author{Su Gao, Arnold W. Miller, and William A. R. Weiss}

\address{Department of Mathematics, PO Box 311430, University of 
North Texas, Denton, TX 76203}
\email{sgao@unt.edu}

\address{Department of Mathematics, University of Wisconsin-Madison,
480 Lincoln Drive, Madison, WI 53706}
\email{miller@math.wisc.edu \hfill www.math.wisc.edu/$\sim$miller}

\address{Department of Mathematics, University of Toronto, Ontario,
CA, M5S 3G3}
\email{weiss@math.toronto.edu}

\begin{abstract}
We prove that there does not exist a subset of the plane $S$ that
meets every isometric copy of the vertices of the unit square in exactly one
point.  We give a complete characterization of all three point subsets $F$ of
the reals such that there does not 
exists a set of reals $S$ which meets every isometric
copy of $F$ in exactly one point.
\end{abstract}



\maketitle

\section{Introduction}

A finite set $X\sq \reals^2$ is {\it Jackson} iff for every
$S\sq\reals^2$ there exists an isometric copy $Y$ of $X$ such that
$|Y\cap S|\neq 1$.  

\begin{question}[Jackson] Is every finite set $X\sq\reals^2$ of two
or more points Jackson? 
\end{question}

This question is motivated by the solution of the Steinhaus problem due to
Jackson and Mauldin \cite{jm1, jm2, jm3}.   They showed that there
exists $S\sq\reals^2$ such that $S$ contains exactly one point from each
isometric copy of $\zz^2$, i.e., there is a Steinhaus set for $\zz^2$. 
Analogous results were obtained by Komjath \cite{kom1,kom2} and
Schmerl \cite{schmerl} for $\zz$, $\qq$, and $\qq^2$.

Note that a set is Jackson iff the Steinhaus problem for it has a negative
solution.  As far as we know the answer to this question is yes, but we have
only partial results.  We do not know if every four point set is Jackson.

In fact we would like to consider a more general version of the problem as the
dimension of the ambient space varies. We will use the following terminology.

\begin{definition} Let $n\geq 1$ and $X\sq\reals^n$. A set $S\sq\reals^n$ is a
{\it Steinhaus set for $X$ in $\reals^n$} if for every isometric copy $Y$ of
$X$ in $\reals^n$, $|Y\cap S|=1$. $X$ is a {\it Jackson set in $\reals^n$} if
there is no Steinhaus set for $X$ in $\reals^n$. 
\end{definition}

In section 2 we focus on finite sets in $\reals^1$ and in section 3 on finite
sets in $\reals^2$. 

The empty set is vacuously Jackson.  A singleton is
never Jackson.
If  $X=\{ x_0, x_1\}\sq \reals$ and $d=|x_0-x_1|>0$, then the set 
 $$S=\displaystyle\bigcup_{n\in\zz}[2nd, (2n+1)d)$$  
is easily seen to be a Steinhaus set for $X$ in $\reals$. Moreover in this case
the Steinhaus set is not unique.  In contrast, if $X=\{ x_0,x_1\}\sq \reals^n$
for any $n\geq 2$, then the following argument shows that $X$ is Jackson.
Suppose $S$ is any Steinhaus set for $X$ in $\reals^n$ and suppose $y_0$ is any
element of $S$. We still use $d(x_0,x_1)$ to denote the distance between $x_0$
and $x_1$. Consider any line $l$ through $y_0$ in $\reals^n$ and let $y_1$ and
$y_2$ be points on $l$ so that $d(y_0,y_1)=d(y_1,y_2)=d$ and $d(y_0,y_1)=2d$.
Since $y_0\in S$ and $\{ y_0,y_1\}$ is an isometric copy of $X$, we have that
$y_1\notin S$. But $\{y_1,y_2\}$ is also an isometric copy of $X$, therefore
$y_2\in S$. Thus we have shown that any point of distance $2d$ from $y_0$ must
be in $S$. In other words, the sphere of radius $2d$ centered at $y_0$ is a
subset of $S$. However, this is a contradiction since it is easy to find two
points on this sphere with distance $d$ between them.  A similar argument shows
that any three point set in $\reals^n$ for $n\geq 2$ is Jackson (see
Proposition \ref{jack1}.) 

This reflection argument is used to derive contradictions in most of
our proofs.

For the rest of the paper we may assume $|X|\geq 3$.

It is easy to see that if a set $X\sq \reals^n$ is Jackson in $\reals^n$ then
as a subset of  $\reals^{n+1}$, $X$ is also Jackson in $\reals^{n+1}$.
The converse is not true.

Another basic point is that if $X$ and $X^\prime$ (in some $\reals^n$) are
similar, then $X$ is Jackson in $\reals^n$ iff $X^\prime$ is Jackson in
$\reals^n$.  This is because we can apply the similarity transformation to a
Steinhaus set for $X$ to get a Steinhaus set for $X^\prime$ (and vice versa).

We end this introduction by connecting our problem with the coloring number
of a distance graph.  For $D$ a set of positive reals define
the graph $G(\reals^n,D)$ by letting $\reals^n$ be the vertices and
connecting $p,q$ by an edge iff $d(p,q)\in D$.  
The chromatic number of a graph $G$ is the smallest $n$ such that
the vertices of $G$ can be partitioned into $n$ sets such that in
each set no two vertices are adjacent.

Given any finite set $F\sq\reals$ define
$$D(F)=\{d(p,q)\;:\; p,q\in F, \;\; p\neq q\}$$

\begin{prop}\label{chrom}
For $F\sq\reals^n$ finite the following are equivalent:
 \begin{enumerate}
  \item There is a Steinhaus set $S\sq\reals^n$ for $F$.
  \item The chromatic number of $G(\reals^n,D(F))\leq |F|$.
 \end{enumerate}
\end{prop}
\proof
Suppose $S$ is a Steinhaus set for $F$.  Consider
$$\{p+S\;:\; p\in F\}.$$
Note that this is a partition of $\reals^n$ into Steinhaus sets for
$F$.  To see that it covers $\reals^n$
note that for any $q\in\reals^n$ that
$(q-F)\cap S\neq \emptyset$
hence there exists $p\in F$ with $q-p\in S$ and so $q\in p+S$.
The sets are pairwise disjoint since otherwise there would
be distinct $p_1,p_2\in F$ and $s_1,s_2\in S$ with $p_1+s_1=p_2+s_2$ 
and so $p_1-p_2=s_2-s_1$ and
so $d(p_1,p_2)=d(s_1,s_2)$.  But
then $S$ would meet an isometric copy of $F$ in at least two points.

Hence the chromatic number of $G$ is less than or equal to
$|F|$.  

Conversely, suppose that the chromatic number of $G$ is $|F|$.  (It cannot be
smaller since the elements of $F$ must
receive different colors.)  Let
 $$\reals^n=\bigcup_{i<|F|}S_i$$
Fix an isometric copy of $F$ say $F'$.
Then each $S_i$ must meet $F'$ in 
at most one point.  But since there are exactly $|F|$ of
the $S_i$ each must meet $F'$ in exactly one point.
\qed

A well-known open problem (see Klee and Wagon
\cite{kw}) is  to determine the chromatic number of the distance one graph in
the plane.  In our terminology this graph would be $G=G(\reals^2,\{1\})$, i.e.,
two points in the plane are adjacent iff the distance between them is exactly
one.  It is known that the chromatic number of $G$ is between 4 and 7.

\begin{prop}
Suppose that the chromatic number of the distance one graph $G$ is 
strictly greater than 4.  Then every four point subset of $\reals^2$ is Jackson.
\end{prop}
\proof
By considering a similar copy of $F$ we may assume there are two points
of $F$ which are exactly one unit apart.  Hence by the above proposition
a Steinhaus set for $F$ would yield a coloring for $G$ with four colors.
\qed

A similar result holds if the chromatic number is greater than 5 or
greater than 6.  Falconer \cite{falconer} has shown that
if the plane is covered by four measurable sets, then one of
the sets contains two points exactly unit one apart.  Hence, no
four point subset of the plane has a measurable Steinhaus set.

\section{Steinhaus sets and Jackson sets on the Line}

All Steinhaus sets and Jackson sets in this section 
will be subsets of $\reals$. 

We begin by giving a complete classification of 3-point Jackson sets on the
line. This is done in two cases. 
By our similarity argument we may assume
that the three point set has the form
$\{0,\alpha,1\}$ where
$0<\alpha<1$.  We begin by considering the
case where $\alpha$ is a rational number. 
In this case (again by similarity) we may
assume that the elements of our set are integers.
Obviously, if $F$ is a finite set of integers, then
$F$ has a Steinhaus set in $\reals$ iff $F$ has a Steinhaus set
$S\sq\zz$, i.e., $S$ meets every isometric copy of $F$ in $\zz$.

\begin{prop} 
\label{rational}
Let $F=\{0,\alpha,1\}$ where $\alpha={n \over m}$
and $n$ and $m$ are relatively prime.
Then the following are equivalent:
\begin{enumerate}
 \item[(i)] There exists a Steinhaus set $S$ for $F$ in $\reals$.
 \item[(ii)] $\{n \mod 3,\;\;\; m \mod 3\}=\{1,2\}$.
 \item[(iii)] $S=\{3k\;:\;k\in\zz\}$ is a Steinhaus set for 
 $\{0,n,m\}$ in $\zz$.
\end{enumerate}
\end{prop}

\begin{proof} 

(ii)$\Rightarrow$(iii):
It is easy to see that every isometric copy of $\{0,n,m\}$ in $\zz$ contains 
exactly one element $u$ such that $u = 0 \mod 3$.

(iii)$\Rightarrow$(i):  Let 
$$S^*=\{u+{k\over m}\;:\;k\in S,\;\;\; 0\leq u < 1\}.$$  
Then since $S$
meets every isometric copy of $\{0,n,m\}$ in $\zz$ in exactly one point,
then $S^*$
meets every isometric copy of $F$ in $\reals$ in exactly one point.

(i)$\Rightarrow$(ii): 
There exists a Steinhaus set $S\sq\zz$ for $\{0,n,m\}$.  And
since the translation of any Steinhaus set is Steinhaus we may assume
$0\in S$.
Note that for any set $S$ its possible periods
$$\{p\in\reals\;:\;\forall x\in\reals\;\; x\in S\rmiff p+x\in S\}$$
is a subgroup of the additive group of $\reals$.  The gaps in the
set $\{0,n,m\}$ are
$a=n$ and $b=m-n$. 
So by the usual reflection argument
(see figure \ref{figreflect}) both $a+2b$ and $2a+b$ are
periods of $S$.  
Let  
$$H=\{ k(2a+b)+l(a+2b)\;:\; k,l\in\zz\}.$$ 
Since $0\in S$ we have that $H\sq S$.
Let $d=gcd(2a+b,a+2b)$ and note that $H=\{kd\;:\; k\in\zz\}$. 
Since $0\in S$ we have that $a\notin S$ and $b\notin S$
and so $d$ does not divide either $a$ or $b$.
 
Note that 
$$a-b=(2a+b)-(a+2b),$$
$$3a=2(2a+b)-(a+2b), \rmand 3b=2(a+2b)-(2a+b),$$
so we have that $d|(a-b)$, $d|3a$ and
$d|3b$. Since
$d$ does not divide $a$, we can write $d=3c$ where $c|a$. 
Similarly,
we have $c|b$. Thus $c=1=\gcd(a,b)$ and $d=3$. 
Now since $d|(a-b)$ we have
$$a\equiv b\,(\mbox{mod}\,3).$$ 
Hence $n\equiv m-n\,(\mbox{mod}\,3)$ but since $n$ and $m$ are relatively 
prime one must be $1\mod 3$ and the other $2\mod 3$.

\begin{figure}
\unitlength=1.36mm
\begin{picture}(90, 40)
\put(20,20){\line(1,0){40}} 
\put(20,20){\circle*{1}} 
\put(30,20){\circle*{1}}
\put(50,20){\circle*{1}} 
\put(60,20){\circle*{1}}
\put(25,15){\text{$b$}}
\put(40,15){\text{$a$}}
\put(55,15){\text{$b$}}
\put(20,10){\vector(0,1){8}}
\put(20,5){\text{$x$}}
\put(60,10){\vector(0,1){8}}
\put(60,5){\text{$x+a+2b$}}
\put(20,40){\line(1,0){50}}
\put(20,40){\circle*{1}} 
\put(40,40){\circle*{1}}
\put(50,40){\circle*{1}} 
\put(70,40){\circle*{1}}
\put(30,35){\text{$a$}}
\put(45,35){\text{$b$}}
\put(60,35){\text{$a$}}
\put(20,30){\vector(0,1){8}}
\put(20,25){\text{$x$}}
\put(70,30){\vector(0,1){8}}
\put(70,25){\text{$x+2a+b$}}
\end{picture}
\caption{Reflection \label{figreflect} }
\end{figure}
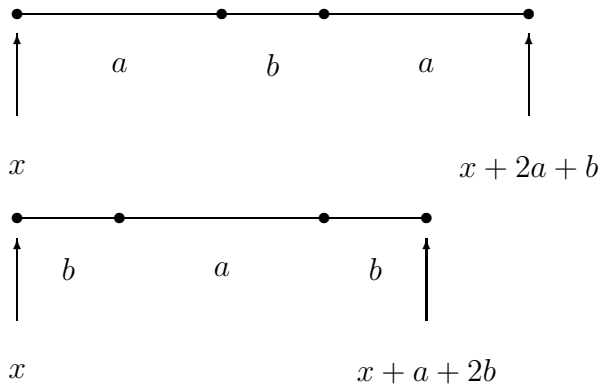

\end{proof}

Another proof of this proposition can be given by using a result
of Zhu \cite{zhu} and the method of Proposition \ref{chrom}. 
Zhu determines the chromatic number of all distance graphs
$G(\zz,D)$ where $|D|=3$.  Our result corresponds to the special case
that $D=\{n,m,m-n\}$  

Note that in Proposition \ref{rational} 
if there is a Steinhaus set for $\{0,{m\over n},1\}$, there is one which
is the union of a closed set and an open set.

\begin{prop} \label{irrational} Let $0<\alpha<1$ be irrational. Then there is a
Steinhaus set for $F=\{0,\alpha,1\}$.  However, there is no Steinhaus set for
$F$ which is either Lebesgue measurable or has the Baire property.
\end{prop}

\begin{proof} 
Let $S$ be a Steinhaus set for $F$.  The gap lengths in $F$ are
$a=\alpha$ and $b=1-\alpha$ and so by the usual reflection argument
(see figure \ref{figreflect}) both $2a+b=1+\alpha$ and $a+2b=2-\alpha$ 
are periods of $S$.  Let  
$$H=\{ m(1+\alpha)+n(2-\alpha)\;:\; m,n\in\zz\}.$$ 
Then for every $h\in H$ and $x\in\reals$ we have that 
$$x+h\in S \rmiff x\in S.$$

Now consider the quotient group $\reals/H$. For any $r\in\reals$, we denote by
$\langle r\rangle$ the coset $r+H$. We claim that $G=\{\langle 0\rangle,
\langle \alpha\rangle, \langle 1\rangle\}$ is a subgroup of $\reals/H$. 
For this note that $1+\alpha\in H$ gives that $\langle
1\rangle+\langle\alpha\rangle=\langle0\rangle$, so
$\langle\alpha\rangle=-\langle1\rangle$. On the other hand, $2-\alpha\in H$
implies that $\langle 2\rangle-\langle\alpha\rangle=\langle0\rangle$, thus
$\langle\alpha\rangle=2\langle 1\rangle$. Putting these together we have that
$2\langle 1\rangle=\langle \alpha\rangle$ and $3\langle 1\rangle=\langle
0\rangle$. Thus the $G$ is a cyclic group of order at most 3. To see that $G$
has order exactly 3, we need only argue that
$\langle1\rangle\neq\langle0\rangle$, i.e., $1\not\in H$. But this is obvious
since otherwise $\alpha$ would be rational.

We are now ready to define a Steinhaus set for $F$.  First
choose a transversal for the cosets of $G$ in $\reals/H$, i.e., a set
$\tilde{S}\subseteq \reals/H$ such that $\tilde{S}$ meets each coset of $G$ at
exactly one point. Then let $S=\bigcup \tilde{S}$, i.e., the union of all
elements of $\tilde{S}$. We check that $S$ is a Steinhaus set as required. For
this let $x\in \reals$. First consider the 3-point set $\{ x, x+\alpha, x+1\}$.
Since the  set $\{\langle x\rangle, \langle x\rangle+\langle \alpha\rangle,
\langle x\rangle+\langle 1\rangle\}=\langle x\rangle+G$ is a coset of $G$ in
$\reals/H$,  there are exactly one element of $\langle x\rangle+G$ which is in
$\tilde{S}$ and therefore exactly one of $x$, $x+\alpha$ and $x+1$ which is in
$S$. Next consider the 3-point set $\{x, x-\alpha,x-1\}$. Using the facts that
$\langle x\rangle-\langle\alpha\rangle=\langle x\rangle+\langle 1\rangle$ and
$\langle x\rangle-\langle1\rangle=\langle x\rangle+\langle \alpha\rangle$ we
can argue similarly that exactly one of $x$, $x-\alpha$ and $x-1$ is in $S$.

Next we show that no Steinhaus set for $F$ can be measurable or have the
property of Baire.
We claim that $H$ is dense in
$\reals$. To see this, 
let $\beta=(1+\alpha)/(2-\alpha)$.  Since $\alpha=(1+\beta)/(2-\beta)$
we have that $\beta$ is
also irrational. 
Let 
 $$K=\{n+m\beta\;:\:n,m\in\zz\}.$$
We claim that $K$ cannot be a discrete subgroup of the reals. If it
were then there would be
some $\epsilon>0$ such that for
every $x,y\in K$ if $|x-y|<\epsilon$
then $x=y$. So if it were discrete there would have to be
minimal positive element $\delta\in K$.
It is easy to see that
 $$K=\{n\delta\;:\:n\in\zz\}$$
since otherwise if $x\in K$ and $n\in\zz$ with
$n\delta<x<(n+1)\delta$ would give us $0<x-n\delta<\delta$.
But since $1\in K$ we have $\delta={1\over n}$ for some positive
integer $n$. Since $\beta\in K$ we would get that $\beta$ is
a rational.  Hence $K$ is dense and since
 $$H=\{(2-\alpha)x\;:\;x\in K\}$$ 
we also have that $H$ is dense.

Since $H+S=S$ and $H$ is dense, it follows that if $S$ has
the property of Baire it must be either meager or comeager.  
This is because if $S$ is comeager in an open interval $I$ then
$h+S$ is comeager in $h+I$.  But $h+S=S$ for any $h\in H$ and
since $H$ is dense, $S$ would have to comeager in $\reals$. 
Similarly, if $S$ is measurable then it is either measure zero
or the compliment of a measure zero set. 

But $S,\alpha+S,1+S$ is a partition of the reals (see the
proof of Proposition \ref{chrom}).  Hence it cannot be either
meager or comeager (or measure zero or comeasure zero).

\end{proof}

As the proofs above show,  the length of the gaps between
successive elements of $F$ are important.

It will sometimes be convenient to use the following gap terminology.

\begin{definition} Let $n\geq 2$ and $a_1,\dots, a_n$ be positive real numbers.
A finite set $X$ of $n+1$ elements {\it has type} $(a_1,\dots, a_n)$ if $X$ is
similar to the set 
$$ A(a_1,\dots,a_n)=\{ 0, a_1, a_1+a_2, \dots, a_1+a_2+\dots+a_n\}. $$ 
\end{definition} 

As remarked before, if $X$ has type $(a_1,\dots,a_n)$, then Steinhaus sets for
$X$ (if they exist) are in one-one correspondence with Steinhaus sets for the
set $A(a_1,\dots,a_n)$ above.  For simplicity we will call the latter sets 
{\it Steinhaus sets for} $(a_1,\dots, a_n)$.

Next we show that Proposition \ref{irrational} can be generalized to
arbitrary sequences $(r_1,\dots,r_n)$ with positive real numbers
$r_1,\dots,r_n$ linearly independent over $\qq$. Note that the linear
independence can be equivalently formulated as follows: for any
$k_1,\dots,k_n\in\zz$, if $k_1r_1+\dots+k_nr_n=0$ then $k_1=\dots=k_n=0$.

\begin{prop} Let $(r_1,\dots,r_n)$ be a sequence of positive real numbers
linearly independent over $\qq$. Then there is a Steinhaus set for
$(r_1,\dots,r_n)$. \end{prop}

\begin{proof} Let $H$ be the additive subgroup of $\reals$ generated by the
elements: 
$$ (n+1)r_1,\;\; r_2-r_1,\;\; r_3-r_1,\;\dots, \;\;r_n-r_1. $$ 
Note 
that these generators are also linearly independent over $\qq$.

Consider the quotient group $\reals/H$ and let $\langle r\rangle$ denote the
coset $r+H$ for $r\in\reals$. Then we note that $\langle r_1\rangle=\langle
r_2\rangle=\dots=\langle r_n\rangle$ and $(n+1)\langle r_1\rangle=\langle
0\rangle$. Let $G$ be the cyclic group in $\reals/H$ generated by $\langle
r_1\rangle$. We claim that $G$ has order $n+1$. To see this, we show that for
any $k\in \zz$, $kr_1\in H$ iff $(n+1)|k$. Let $k_1,\dots, k_n\in \zz$. Then

$$\begin{array}{rl} & kr_1=k_1(n+1)r_1+k_2(r_2-r_1)+\dots+k_n(r_n-r_1) \\
\Leftrightarrow & (k_1(n+1)-k-k_2-\dots-k_n)r_1+k_2r_2+\dots+k_nr_n=0 \\
\Leftrightarrow & k_1(n+1)-k-k_2-\dots-k_n=k_2=\dots=k_n=0 \\
\Rightarrow & k_1(n+1)=k.
\end{array}
$$

We construct $S$ as before. First let $\tilde{S}$ be a transversal for the
cosets of $G$ in $\reals/H$. Then let $S=\bigcup\tilde{S}$. Then $S$ is a
required Steinhaus set by a similar argument as before. 
\end{proof}

The following example shows that there are Jackson sets of type
$(r_1,\dots,r_n)$ where at least one pair of numbers have an irrational ratio.

\begin{example} For any real number $r>0$, the set $A(1,r,2)$ is Jackson. 
\end{example}

\begin{proof} Assume $S$ is a Steinhaus set for $A(1,r,2)=\{0,1,1+r,3+r\}$ 
and $0\in S$. 
Then $-1, 1, 1+r, 2+r$ are each not in $S$ since their distances to $0$  
are forbidden. i.e, in $D(A(1,r,2))$.  But
the set $\{-1,1,1+r,2+r\}$ is an isometric copy of $\{0,1,1+r,3+r\}$, 
(i.e. reflect across the interval of length $r$).  Thus it
should meet $S$ in exactly one point, a contradiction. 
\end{proof}

In what follows we will give various other examples of 4-point Jackson sets. We
have ad hoc arguments for  4-point sets of more types than covered here.
However, a complete classification for 4-point Jackson sets is not known. 

For the remainder of this section we will only consider finite subsets of the
integers $\zz$ and Steinhaus sets in $\zz$.

Recall that
a set $S\sq \reals$ is {\it periodic} if there is $p>0$ such that 
for any $x\in\reals$, 
 $$x\in S \rmiff x+p\in S.$$ 
In this case $p$ is called a {\it period} for $S$.

\begin{prop} Let $F\sq\zz$ be finite and let 
$d$ be the diameter of $F$, i.e., $d=\max(F)-\min(F)$.
If $F$ has a Steinhaus set, then $F$ has a Steinhaus set
with integer period $p$
with $0<p\leq d2^d$.  
\end{prop}

\begin{proof}

Let $S$ be a Steinhaus set for $F$. For $0\leq k\leq 2^d$, define
a function $f_k:\{0,\dots,d-1\}\to\{0,1\}$ by letting 

$$f_k(t)=1\ \mbox{ iff }\ ks+t\in S.$$ 

\noindent Since there are $2^d+1$ functions $f_k$ but only $2^d$ many functions
from $\{0,\dots,d-1\}$ to $\{0,1\}$, there must be $k_1<k_2$ such that
$f_{k_1}=f_{k_2}$. 

Let $p=(k_2-k_1)d$ and define a set $S'\sq\zz$ to
be periodic with period $p$ such that
$S'\cap I=S\cap I$ where $I$ is the interval $[k_1d,(k_2+1)d)$.
This is possible because $S$ is the same on the first and last $d$-subintervals
of $I$.

We check that $S'$ is a Steinhaus set for $F$. Let $X$ be an isometric copy of
$F$ in $\zz$. 
Let $x_0$ be the smallest element of $X$ and note that $x_0+d$ is the
largest element of $X$.  Choose $n\in\zz$ so that
  $$(n-1)p+x_0<k_1d\leq np+x_0.$$ 
Note that $np+x_0<k_2d$.  Hence,
$np+X\sq I$.  Since $S$ and $S'$ agree on $I$
it must be that $S'$ meets $np+X$ in exactly one point.  Since
$S'$ has period $p$ it meets $X$ in exactly one point.

\end{proof}

A proof of this proposition can also be given by noting that periodic colorings
give a periodic Steinhaus set (proof of Proposition \ref{chrom}) and using the
proof of Theorem 2 of Eggleton, Erdos, and Skilton \cite{egg}.   Their proof
works for any finite distance set $D\sq \zz$ although they assume $D$ contains
only primes.

Let $\mathcal J$ denote the set of all finite sequences of positive integers
$F$ such that $F$ is Jackson.

\begin{theorem}
The set $\mathcal J$ is computable.
\end{theorem}

\begin{proof} 
It is clear from the last proof that $F$ with
diameter $d$ has a Steinhaus set iff
there exists $S\sq\{0,1,\ldots,d2^d\}$ which meets every isometric copy
of $F$ inside $\{0,1,\ldots,d2^d\}$ in exactly one point.
But this is clearly computable.
\end{proof}

However, the following question seems to be open.

\begin{question} Is $\mathcal J$ computable in polynomial time?
\end{question}

Regarding periodic Steinhaus sets we have the following useful fact.

\begin{prop}\label{period} Let $F$ be a finite set of positive
integers. If $S$ is a periodic Steinhaus set for $F$,
with an integer period $p$, then $p$ is a multiple of $n=|F|$.  In fact,
$p=mn$ where $m=|S\cap\{0,1,\ldots,p-1\}|$.  
\end{prop}

\begin{proof} 
Let $F=\{b_1,b_2,\ldots,b_n\}$.
Consider the $p\times n$-matrix $(c_{ij})$ defined as
follows: for $i=0,\dots, p-1$ and $j=1,\dots, n$, let 
$$c_{ij}=\left\{
\begin{array}{ll}
1 & \mbox{ if } i+b_j\in S\\
0 & \mbox{ otherwise }\\
\end{array}\right.
$$
Note that for each $i<p$ there is exactly one $j$ such that
$c_{ij}=1$ because $|(i+F)\cap S|=1$. Hence the total
number of ones in the matrix $C$ is $p$.

On the other hand for each fixed $j=1,\ldots, n$
$$|\{0\leq i<p\;:\; i+b_j\in S \}|=|\{0\leq i<p\;:\; i\in S \}|=m$$
since $p$ is a period of $S$.
Hence the total number of ones in the matrix $C$ is $mn$. 

Therefor $p=mn$.

\end{proof}

The proposition has many applications. One can get a taste of the flavor of
these applications from the simple example below.

\begin{example}\label{example1} Let $a$ and $b$ be positive integers such that
$a+b$ is odd. Then the set $A(a,b,b)$ is Jackson. 
\end{example}

\begin{proof} Toward a contradiction assume $S$ is a Steinhaus set for
$A(a,b,b)$.   Note that by reflection $2a+2b=a+b+b+a$ is a period for $S$.
But since $a+b$ is odd, $2a+2b$ is not a multiple of $4$, 
contradicting the preceding proposition. 
\end{proof}

In the next few propositions we investigate when a Steinhaus set exists with a
given period.  We start with the simplest case.

\begin{prop}\label{1time} Let $(a_1,\dots,a_n)$ be a sequence of positive
integers. The following are equivalent:
\begin{enumerate}
\item[(i)] A Steinhaus set of period $n+1$ exists for $(a_1,\dots,a_n)$. 
\item[(ii)] $A(a_1,\dots,a_n)\equiv\{0,\dots,n\}\ (\mbox{\rm mod}\,n+1)$.
\item[(iii)] For any $x\neq y$ elements of $A(a_1,\dots,a_n)$, $x\not\equiv y\
(\mbox{\rm mod}\,n+1)$.
\end{enumerate}
\end{prop}

\begin{proof} Clauses (ii) and (iii) are obviously equivalent since
 $$|A(a_1,\dots,a_n)|=|\{0,\dots,n\}|=n+1$$ 
and thus a map (in this context the
mod $n+1$ map) between the two sets is onto iff it is one-one. To see that
(i)$\Rightarrow$(iii), let $S$ is a Steinhaus set of period $n+1$ for
$(a_1,\dots,a_n)$ and without loss of generality assume $0\in S$.  Assume that
there are distinct $x$ and $y$ from $A(a_1,\dots,a_n)$ with $x\equiv y\
(\mbox{mod}\,n+1)$. Then $y-x\in S$ by the periodicity of $S$.  Now $0, y-x\in
A(a_1,\dots,a_n)-x$ and hence $A(a_1,\dots,a_n)-x$, an isometric copy of
$A(a_1,\dots,a_n)$, meets $S$ at two points, a contradiction.

To see that (ii)$\Rightarrow$(i), it suffices to note that the set 
$$S=\{k(n+1)\;:\; k\in\zz\}$$
is a Steinhaus set for
$(a_1,\dots,a_n)$. $S$ obviously has period $n+1$. Let $X\sq\zz$ be 
any isometric
copy of $A(a_1,\dots,a_n)$. Then the mod $n+1$ values of 
$X$ are distinct. It follows that there is exactly one element 
of $X$ in $S$.
\end{proof}

\begin{example} There are Steinhaus sets of period $4$ for $(a,b,c)$ iff
$(a,b,c)$ is congruent mod $4$ to one of the following triples: 
$$ (1,1,1),\ (1,2,3),\
(2,3,2),\ (2,1,2),\ (3,2,1),\mbox{ and }(3,3,3).$$ 
\end{example}

\begin{proof} This follows from the preceding proposition by a direct
computation. 
\end{proof}

\begin{definition}
Given a sequence $(a_1,\dots,a_n)$ $(n\geq 2)$ of positive integers.
Define 
$$\dist=\dist(a_1,\dots,a_n)=\{ x-y\,|\, x,y\in A(a_1,\dots,a_n)\}.$$ 
\end{definition}

It is convenient to include in $\dist$ 
the negatives of the forbidden distances. 
Note that if $S$ is a Steinhaus set for $(a_1,\dots,a_n)$ and $0\in S$
(this does not lose generality since any shift of $S$ is still a Steinhaus
set), then $\dist\cap S=\{0\}$. 
For a putative period $M$ we denote by $\dist\,(\mbox{mod}\,M)$ the set of mod
$M$ values of elements of $\dist$.

\begin{prop}\label{pSM} 
Let $(a_1,\dots,a_n)$ be a sequence of positive integers and let $k\geq 0$.
Then there exists a Steinhaus set for
$(a_1,\dots,a_n)$ of period $M=(k+1)(n+1)$  iff
\begin{enumerate}
\item[(a)] Elements of $A(a_1,\dots,a_n)$ have distinct mod $M$ values; and
\item[(b)] There are integers, $0=x_0<x_1<\dots<x_k <M$, such that 
$$x_j-x_i\notin \dist\,(\mbox{\rm mod}\,M)$$ 
for all $i<j\leq k$.
\end{enumerate}
\end{prop}

\begin{proof} First assume that $S$ is a Steinhaus set for
$(a_1,\dots,a_n)$ of period $M$  and without loss of generality $0\in S$. 
Clause (a) follows
from a similar argument as in the preceding proof.  It remains to show (b).
From Proposition \ref{period} we know that
$|S\cap\{0,\dots,M-1\}|=k+1$. 
Let $S\cap\{0,\dots, M-1\}=\{0, x_1,\dots,x_k\}$.
Suppose for contradiction that $x_j-x_i=u-v+mM$ for some
$m\in\zz$ and distinct $u,v\in A(a_1,\dots,a_n)$.  But by
the periodicity of $S$ we have that $x_j-mM\in S$ and hence
we have two elements $x_i, x_j-mM$ of $S$
at a forbidden distance $d(u,v)$. 

Conversely, assume that (a) and (b) hold. Let $S$ be periodic with period $M$
and $S\cap [0,M)=\{0,x_1,\ldots,x_k\}$ 
We check that $S$ is a Steinhaus set for $A(a_1,\dots,a_n)$
with respect to $\zz$.  Let $X$ be
an isometric copy of $A(a_1,\dots,a_n)$.
First we argue that $|X\cap S|\leq 1$. For this let
$x\neq y\in X$ be both in $S$. 
Let $x'=x\,(\mbox{mod}\,M)$ and $y'=y\,(\mbox{mod}\,M)$. Then by (a) $x'\neq
y'$ and  by periodicity of $S$ we have that $x',y'\in S$. So let $x'=x_j$ and
$y'=x_i$.  But then,  $x_j-x_i\in \dist\,(\mbox{\rm mod}\,M)$ contradicting
(b).

To see that $X\cap S\neq\emptyset$, we first consider the case that $X$ is a
shift of $A(a_1,\dots,a_n)$, i.e., there is some $c\in\zz$ such that
$X=c+A(a_1,\dots,a_n)$. 
Since $M$ is a period of $S$, 
for any $b\in\zz$ we have that 
$$|(b+S)\cap [0,M)|=k+1.$$  
Since
$S$ meets each isometric copy of $A(a_1,\dots,a_n)$ in at most
one point we have that 
$$\{-b+S\;:\; b\in A(a_1,\dots,a_n)\}$$ 
are
pairwise disjoint.  Since $M=(k+1)(n+1)$ and 
$$|A(a_1,\dots,a_n)|=n+1$$
it follows that $\{-b+S\;:\; b\in A(a_1,\dots,a_n)\}$ must cover
$[0,M)$ and hence all of $\zz$. But if $c\in -b+S$
then 
$$c+b\in S\cap (c+A(a_1,\dots,a_n)).$$ 
A similar argument can be given for the case that $X=c-A(a_1,\dots,a_n)$.
\end{proof}

\begin{corollary}\label{2times} Let $(a_1,\dots,a_n)$ be a sequence of positive
integers. Then there exists a Steinhaus set of period $2(n+1)$ for
$(a_1,\dots,a_n)$ iff

\begin{enumerate}
\item[(a)] Elements of $A(a_1,\dots,a_n)$ have distinct mod $2(n+1)$ values; 
and
\item[(b)] $\dist\not\equiv\{0,\dots,2n+1\}\ (\mbox{\rm mod}\,2(n+1))$.
\end{enumerate}
\end{corollary}

\begin{proof} Clause (b) is equivalent to saying that there is an $x$ such
that $x\not\in\dist\,(\mbox{mod}\,2(n+1))$. 
\end{proof}

Let $k\geq 0$ be an integer.  Sets of type $(1,1,4k+2)$ or $(1,1,4k+4)$
are Jackson by Example \ref{example1}.  Sets of type $(1,1,4k+1)$ have Steinhaus
sets of period 4 by Proposition \ref{1time}.

\begin{example} Any set of type $(1,1,4k+3)$ is
Jackson.  
\end{example}

\begin{proof}  
The reflection argument gives that if $S$ is a Steinhaus set for 
$(1,1,4k+3)$ then
$S$ has a period $8(k+1)$. When $k=0$ the above corollary applies. But in this
case $\dist=\{0,1,2,3,4,5,-1,-2,-3,-4,-5\}$ and hence
$\dist\,(\mbox{mod}\,8)=\{0,1,2,3,4,5,6,7\}$. Thus there does not exist any
Steinhaus set of period $8$ for $(1,1,3)$ and therefore any set of type
$(1,1,3)$ is Jackson.

For the general case $k>0$ we let $M=8(k+1)$ and get that 
$$\dist=\{ 0,\pm 1,\pm 2,\pm (4k+3),\pm (4k+4), \pm (4k+5)\},$$ 
and therefore 
$$\dist\,(\mbox{mod}\,M)=\{ 0,1,2, 4k+3,4k+4,4k+5, 8k+6,8k+7\}. $$ 
Write
$B=[3,4k+2]\cap\nn$ and $C=[4k+6,8k+5]\cap\nn$. Then $B\cup C$ is the
complement of $\dist\,(\mbox{mod}\,M)$ in $\{0,\dots,M-1\}$. Assume that a
Steinhaus set of period $M$ for $(1,1,4k+3)$ existed. Then by Proposition
\ref{pSM} there are $x_1,\dots,x_{2k+1}\in B\cup C$ with $x_i-x_j\in B\cup
C\,(\mbox{mod}\,M)$. Let $\{x_1,\dots,x_{2k+1}\}\cap B=B_0$, $h=|B_0|$,
$C_0=\{x_1,\dots,x_{2k+1}\}\cap C$ and $l=|C_0|$. Then $h+l=2k+1$. We note the
following facts. For each $x\in B_0$ such that $x\leq 4k$, the numbers $4k+3+x,
4k+4+x$ and $4k+5+x$  are in $C\cap (x+\dist)$, hence in particular are not in
$C_0$. Moreover, since $1,2\in \dist$, any two elements of  $B_0$ differ by at
least 2, the two sets of three numbers associated with them have an empty
intersection.  Using this we obtain the following inequality $$ l+(h-1)\cdot
3+1\leq 4k, $$ where the left hand side counts the number of elements in $C_0$
and the number of elements excluded from $C_0$ for being associated with one of
the numbers in $B_0$ (in a worst case scenario), and the right hand side the
size of $C$. By a symmetric and similar argument we also obtain $$ h+(l-1)\cdot
3+1\leq 4k,$$ by reversing the roles of $B$ and $C$. Thus we get $$
4(h+l)-4\leq 8k \mbox{ or } h+l\leq 2k+1. $$ This is not a contradiction yet
but it follows from $h+l=2k+1$ that all elements of $B$ and $C$ are accounted
for in the above counting in order to establish the inequalities. This is to
say that every element of $B$ is either an element of $B_0$ or else is
associated with an element of $C_0$, and vice versa for elements of $C$. Now we
claim that either $3\in B_0$ or $4k+6\in C_0$. It is easy to see that if
$4k+6\not\in C_0$, then the only element of $B$ with which it is associated is
$3$. In the case $3\in B_0$ it is straightforward to check that $3+4i\in B_0$
for all $i=1,\dots,k-1$ and $4k+5+4j\in C_0$ for all $j=1,\dots,k$, but this is
a contradiction since there are only $2k$ of them. Similarly in the case
$4k+6\in C_0$ it is straightforward to check that $4i\in B_0$ for $i=1,\dots,k$
and $4k+6+4j\in C_0$ for $j=1,\dots,k-1$, which is a total of $2k$ of them,
contradiction again. \end{proof}

Note that the Proposition \ref{1time} and Corollary \ref{2times} give
polynomial time computable criteria for the existence of Steinhaus sets of
period $n+1$ and $2(n+1)$, respectively. However, the criterion in Proposition
\ref{pSM} is NP in $(a_1,\dots,a_n)$ and $M$. We do not know whether the
existence of Steinhaus sets for a given period is computable in polynomial
time. The following corollary provides an easy sufficient condition for this
existence.

\begin{question}
Determine all pairs of positive integers $a,b$ such that sets of
type $(a,a,b)$ are Jackson.
\end{question}

\begin{corollary} Let $(a_1,\dots,a_n)$ be a sequence of positive integers and
let $s=a_1+\dots+a_n$. If $M=(k+1)(n+1)>s$ and there is an integer $0<x<M$ such
that 
 $$\{x,2x,\dots,kx\}\cap \dist=\emptyset (\mbox{mod}\,M),$$ 
then there is a Steinhaus set of
period $M$ for $(a_1,\dots,a_n)$. 
Consequently, if
$(b_1,\dots,b_n)\equiv(a_1,\dots,a_n)\,(\mbox{\rm mod}\,M)$, with $(a_1,\dots,
a_n)$ and $M$ satisfying the conditions above, then there is also a Steinhaus
set of period $M$ for $(b_1,\dots,b_n)$. 
\end{corollary}

\begin{proof} If $M>s$ then clause (a) of Proposition \ref{pSM} holds. Let
$x_i=ix$ for $i=0,1,\dots,k$. Then 
for any $0\leq i<j\leq k$, 
$$x_j-x_i=(j-i)x=x_{j-i}\notin\dist\,(\mbox{mod}\,M).$$ 
Thus clause (b) of Proposition \ref{pSM} is also
verified. The additional part is trivial.  However, the point is that
$b_1+\dots+b_n$ may not be less than $M$. 
\end{proof}

\begin{example} Let $a$ and $b$ be odd positive integers with
$a+b\equiv0\,(\mbox{\rm mod}\,4)$. Then there is a Steinhaus set of period
$2(a+b)$ for $(a,b,a)$. 
\end{example}

\begin{proof} Let $a+b=4k$ and $M=2(a+b)=8k$. Then
$\dist\,(\mbox{mod }8k)\,=\{0,a,4k-a, 4k, 4k+a,8k-a\}$. We have that
$\{2,4,\dots,4k-2\}\cap \dist=\emptyset$ since both $a$ and $4k-a=b$ are odd. 
The proof is completed with the application of the preceding corollary. 
\end{proof}

\section{Jackson Sets on the Plane}

\begin{prop}[Jackson] \label{jack1}
 Every three element subset of the plane is Jackson.
\end{prop}
\begin{proof}
For contradiction suppose there were such an $S$.
Let $X=\{p,q,r\}$ and suppose $p\in S$.  Form the parallelogram
$pqrp^\prime$.  Now if $p\in S$ then $r,q\notin S$ and therefor 
$p^\prime\in S$ since the triangle
$\triangle pqr$ is isometric to $\triangle p^\prime rq$. 
Let $C$ be the circle centered at $p$
of radius $d(p,p^\prime)$.  It follows that every point on $C$ is in $S$, but
this is a contradiction since there is a chord of $C$ with length
$d(q,r)$ and therefor there would be an isomorphic copy of $X$ containing
at least two points of $S$.
\end{proof}

Note that this propositional also holds for any three point set in $\reals^n$
for $n>2$.  A similar argument reflection argument would work for four
noncollinear points in $\reals^3$ provided that $4a\geq e$  where $a$ is
maximum altitude of the tetrahedron and $e$ the minimum edge length.

\begin{prop}[Jackson]\label{jack2}
 Every set of $n>1$ equally spaced collinear points 
is Jackson.
\end{prop}
\begin{proof}
Suppose the points are $p_0,p_1,\ldots,p_n$ written in their
natural order.  
Let $p_{n+1}$ be the next point along the line containing
$p_0,p_1,\ldots,p_n$ with  the same distance apart.   
Suppose that $p_0\in S$, then $p_1,\dots,p_n\notin S$ and 
since the set $p_1,\ldots,p_{n+1}$ is isometric to the original
it must be that $p_{n+1}\in S$.  It follows that the circle of
radius $d(p_0,p_{n+1})$ around $p_0$ is a subset of $S$ and so we get the
same contradiction as in the last proposition.
\end{proof}

This proposition also holds if all gaps but the first gap are equal.

\begin{figure}
\unitlength=1.36mm
\begin{picture}(90, 40)
\put(10,10){\line(1,1){10}}
\put(10,10){\line(2,1){20}}
\put(30,20){\line(1,1){10}}
\put(20,20){\line(2,1){20}}
\put(10,10){\circle*{1}} \put( 7, 7){\text{$p$}}
\put(20,20){\circle*{1}} \put(17,23){\text{$q$}}
\put(30,20){\circle*{1}} \put(33,17){\text{$r$}}
\put(40,30){\circle*{1}} \put(43,33){\text{$p^\prime$}}
\multiput(20,20)(2,0){5}{\line(1,0){1}}
\put(60,20){\line(1,0){30}} 
\put(60,20){\circle*{1}} \put(60,16){\text{$p_0$}}
\put(65,20){\circle*{1}} \put(65,16){\text{$p_1$}}
\put(70,20){\circle*{1}} \put(70,16){\text{$p_2$}}
\put(75,20){\circle*{1}} \put(78,16){\text{$\ldots$}}
\put(85,20){\circle*{1}} \put(85,16){\text{$p_n$}}
\put(90,20){\circle*{1}} \put(90,16){\text{$p_{n+1}$}}
\end{picture}
\caption{Proposition \ref{jack1} and \ref{jack2}}
\end{figure}

\begin{prop} 
A finite set $X\sq \reals^2$ is Jackson iff there exists a finite set 
$F\sq\reals^2$ such that for every $S\sq F$ there exists $Y\sq F$ an isometric
copy of $X$ such that $|S\cap Y|\neq 1$. 
\end{prop}
\proof
This follows from the compactness theorem. Let 
$$\pp=\{P_z \st z\in\reals^2\}$$
be propositional letters and let
$$\cc=\{Y\sq\reals^2\st Y \mbox{ is an isometric copy of } X\}.$$ 
For each $Y\in\cc$ let $\th_Y,\psi_Y$ be the propositional sentences:
$$\th_Y=\bigvee_{z\in Y}P_z$$
$$\psi_Y=\bigwedge_{z,w\in Y, z\neq w}(\neg P_z\vee \neg P_w)$$
i.e. these say that $P_z$ holds for exactly one $z\in Y$.
Let
$$T=\{\th_Y,\;\psi_Y\st Y\in \cc\}.$$
A truth evaluation (model) for $T$ corresponds exactly to
a subset $S\sq\reals^2$ witnessing that $X$ is not Jackson.
By the compactness theorem for propositional logic $T$ since
$T$ is inconsistent there is a finite subset of $T$ which is inconsistent.
Hence there exist $Y_1,\ldots,Y_n\in\cc$ such that
$$\Sigma=\{\th_{Y_1},\;\psi_{Y_1},\;\th_{Y_2},\;\psi_{Y_2},
\;\ldots,\; \th_{Y_n},\;\psi_{Y_n}\}$$
is inconsistent and therefor $F=\cup_{i=1}^n Y_i$ does the trick.
\qed 

Another proof (but a little round-about) is to quote the
Erdos-Bruijn Theorem \cite{eb} about graph coloring.

\begin{prop}
The vertices of the unit square is Jackson.  More generally, the
vertices of a rectangle with
height $x$ and width 1 where $x^2$ is a rational number is Jackson. 
\end{prop}
\proof
First consider the case of the vertices of the unit
square.  For contradiction assume
that $S\sq\reals^2$ meets every isomorphic
copy the vertices of the unit square in exactly one point.
Consider $S\cap \zz^2$.  There are two possibilities. 

\bigskip\noindent
{\bf Claim 1.} Either
\begin{itemize}
  \item $S$ is disjoint from every other column of $\zz^2$ and for
   the other columns it contains every other point.
  \item $S$ is disjoint from every other row of $\zz^2$ and for
   the other rows it contains every other point.
\end{itemize}
In either case $S\cap\zz$ does not contain two points in the same
row or column and at an odd distance.

\bigskip
To see  why this is so, suppose the point $(i,j)\in S$. Then the
points $(i+1,j),(i,j+1),(i+1,j+1)\notin S$.  But since the square
next to it must meet $S$, exactly one of two points $(i+2,j)$
or $(i+2,j+1)$ must be in $S$.  

\begin{figure}
\unitlength=1mm
\begin{picture}(60, 80)(-30,0)
\multiput(10,10)(0,20){4}{\line(1,0){40}}
\multiput(10,10)(20,0){3}{\line(0,1){60}}
\put(10,10){\circle*{2}}
\put(50,30){\circle*{2}}
\put(10,50){\circle*{2}}
\put(50,70){\circle*{2}}
\put(10,6){\text{$i$}}
\put(30,6){\text{$i+1$}}
\put(50,6){\text{$i+2$}}
\put(2,10){\text{$j$}}
\put(2,30){\text{$j+1$}}
\put(2,50){\text{$j+2$}}
\put(2,70){\text{$j+3$}}
\put(40,40){\text{$Q$}}
\end{picture}
\caption{Knight's move\label{kn}}
\end{figure}
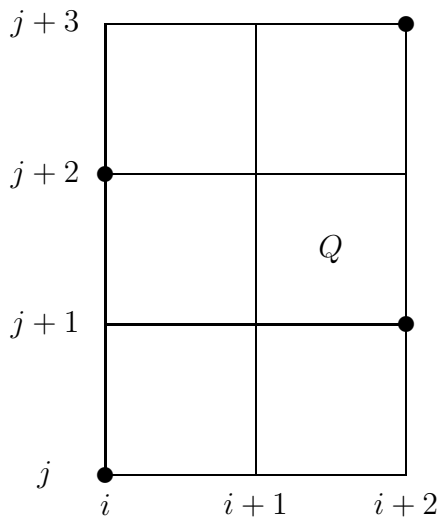

Case 1. Suppose $(i,j)\in S$ and
$(i+2,j+1)\in S$, i.e., the knight's move in chess, see figure
\ref{kn}.  

In this case
the points $(i,j+2)$ and $(i+2,j+3)$ must be in $S$.  
This is because the square $Q$
with left lower corner $(i+1,j+1)$ rules out the point $(i+1,j+2)$
and the point $(i+2,j+2)$.  Continuing with this argument we see
that $S$ must contain all the points $(i,j+2n)$ and $(i+2,j+2n+1)$ for
$n\in\zz$. It must also be disjoint from the columns $i-1,i+1,i+3$.      
An empty column implies that the next column must contain every other
point and a column which contains every other point implies that the
next column is disjoint.  

Case 2. Suppose $(i,j)\in S$ and $(i-2,j+1)\in S$.

This is symmetrical and again we get that every other column is disjoint
from $S$.

Case 3. Whenever $(i,j)\in S$ we have $(i\pm 2, j)\in S$.

In this case every row $\{i+2n+1\}\times\zz$ is disjoint from $S$ and
so on.

This proves Claim 1.  Note that by symmetry it is true for every
lattice (isometric copy of $\zz^2$).

\bigskip\noindent
{\bf Claim 2.}  Suppose $A\in S$ then the circle of radius $4$ around $A$ is
a subset of $S$ and so we are done.
\bigskip
 
Consider any three points $A,B,C$ with $d(A,B)=4$,
$d(B,C)=3$ and $d(A,C)=5$.  To simplify notation assume that
$A=(0,0)$, $B=(0,4)$, and $C=(3,4)$.  

Case 1. For every odd $n\in\zz$ we have $(\{n\}\times\zz)\cap S=\emp$.

This means in particular that the point $D=(1,4)\notin S$.  But
since either $B$ or $D$ is in $S$ we have that $B\in S$.

Case 2. For every odd $n\in \zz$ we have $(\zz\times\{n\})\cap S=\emp$.

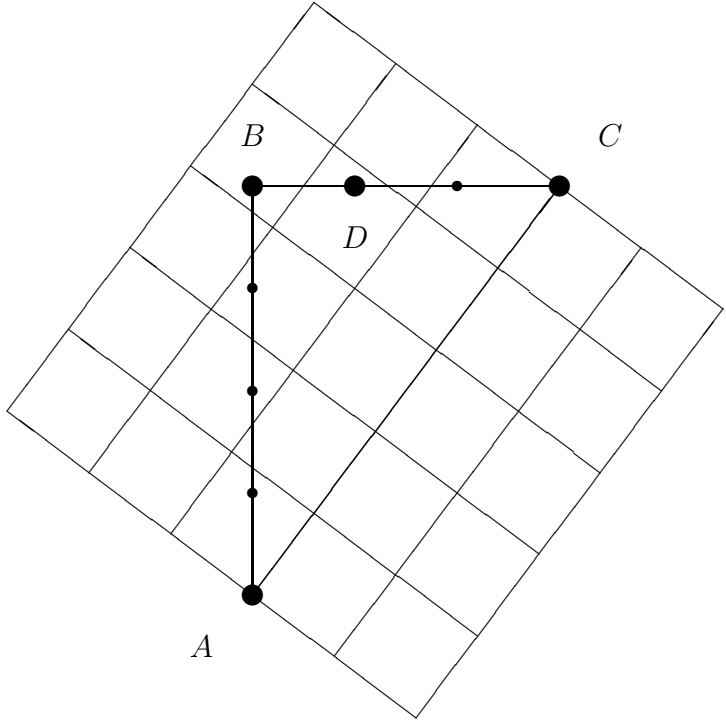
\begin{figure}
\unitlength=1.36mm
\begin{picture}(90,80)(0,15)
\put(30,30){\circle*{2}}    \put(25,25){\text{$A$}} 
\multiput(30,40)(0,10){3}{\circle*{1}}
\put(30,30){\line(0,1){40}}
\put(30,30){\line(3,4){30}}
\put(30,70){\circle*{2}}    \put(30,75){\text{$B$}}
\put(40,70){\circle*{2}}    \put(40,65){\text{$D$}}
\put(50,70){\circle*{1}}
\put(30,70){\line(1,0){30}}
\put(60,70){\circle*{2}}    \put(65,75){\text{$C$}}
\multiput(46,18)(-8,6){6}{\line( 3,4){30}}
\multiput(46,18)( 6,8){6}{\line(-4,3){40}}
\end{picture}
\caption{Tipped Lattice $\ll$\label{tl}}
\end{figure}

In this case every other point on the row $\zz\times\{4\}$ is in $S$.
If $B\notin S$ then the point $D=(1,4)\in S$ and so is the point
$C=(3,4)$.  But consider the tipped lattice $\ll$ (see figure \ref{tl})
with the points
$(0,0)$ and $(3,4)$ where these points are
on the same ``column'' of $\ll$. 
$S$ must satisfy Claim 1 for $\ll$ but
$A$ and $C$ are on the same column of $\ll$ and separated by
an odd distance (5) which is a contradiction.
This does the special case of the unit square.  

\bigskip
Now suppose our rectangle has height $x$ and width $1$
and $x^2$ is a rational number. 

Here the basic lattice to consider is
$$\{(m,nx)\st m,n\in\zz\}.$$

The parity argument we have just given works for any $x$ satisfying
$$(odd)^2+(even)^2x^2= (odd)^2$$
for some positive odd and even integers, say $a^2+b^2x^2=c^2$. 
Consider 3 points $A$, $B$, and $C$ with 
$d(A,B)=bx$, $d(B,C)=a$, and $d(A,C)=c$.  We claim that
if $A\in S$, then $B\in S$.  Hence again we get a contradiction
since the circle of radius $bx$ centered at $A$ is a subset of
$S$.  So suppose for contradiction that $B$ is not in $S$.  
Then since the distance from $A$ to $B$ is $bx$ which is even number
times the height of the rectangle and since $B\notin S$, 
it must be that $S$ contains every
other element of the row on which $B,C$ are on.
But since $a$ is odd number this means 
$C\in S$.  Now consider the
tipped $1,x$-lattice $\ll$ with the $1$ side along
the line segment joining $A$ to $C$.  But this length is $c$ which is
odd and this cannot 
happen in a lattice satisfying the analogue of Claim 1. 

Fix any positive rational number $k/l$.  If 
$$(2n-1)^2+(2m)^2x^2=(2n+1)^2$$
then 
$$x^2={{(2n+1)^2-(2n-1)^2}\over{(2m)^2}}={{8n}\over{(2m)^2}}={k\over l}$$
if we take $n=2kl$ and $m=2l$.
\qed

\begin{question}
Is the set of vertices of a rectangle of width $^3\sqrt{2}$ and height 
$1$ Jackson?
\end{question}

\begin{question}
Does there exists $S\sq\nn$ such that $S$ contains exactly one
element of every Pythagorean triple $x^2+y^2=z^2$?
\end{question}

It might be helpful to note that two Pythagorean triples can have
at most one point in common, i.e., it impossible to have
$$x^2+y^2=z^2 \rmand y^2+z^2=w^2$$
in the positive integers, see Sierpi\'nski \cite{sier}
\S 10.6.

\begin{question}(Steprans) Does there exist a finite partition of
$\nn$, say $\nn=S_1\cup S_2\cup\cdots\cup S_n$ such that no
$S_i$ contains a Pythagorean triple?
\end{question}

\begin{prop}
For any trapezoid $T$ and $\epsilon>0$ there exists a trapezoid
$T^\prime$ with corresponding vertices within $\epsilon$ of the 
vertices of $T$ such that the vertices of $T^\prime$ are Jackson.
\end{prop}
\proof
Consider any trapezoid with sides $A,B,C,D$.  (By trapezoid
we just mean that two of the sides are parallel.)
Assembly four isometric copies into the chevron. 

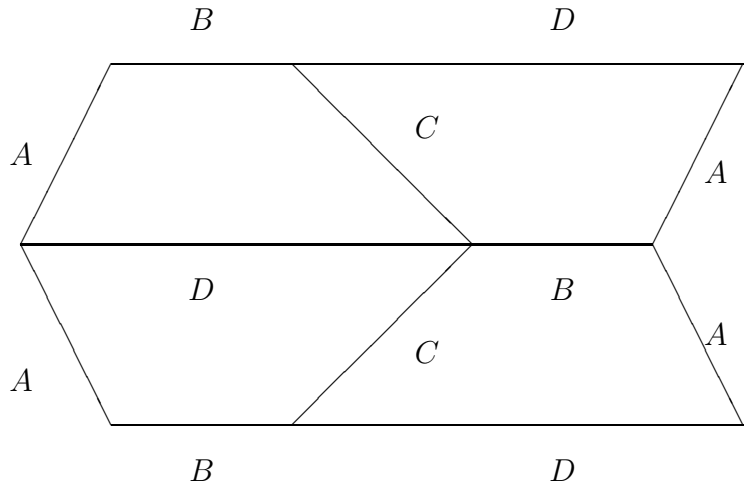
\begin{figure}
\unitlength=1.2mm
\begin{picture}(100, 70)(10,10)
\put(20,40){\line(1,0){70}} 
\put(20,40){\line(1, 2){10}}
\put(20,40){\line(1,-2){10}}
\put(30,60){\line(1, 0){70}} 
\put(30,20){\line(1, 0){70}} 
\put(90,40){\line(1, 2){10}}
\put(90,40){\line(1,-2){10}}
\put(50,60){\line(1,-1){20}}
\put(50,20){\line(1, 1){20}}

\put(20,25){\text{$A$}}
\put(20,50){\text{$A$}}

\put(40,15){\text{$B$}}
\put(40,35){\text{$D$}}
\put(40,65){\text{$B$}}

\put(65,28){\text{$C$}}
\put(65,53){\text{$C$}}

\put(80,15){\text{$D$}}
\put(80,35){\text{$B$}}
\put(80,65){\text{$D$}}

\put(97,30){\text{$A$}}
\put(97,48){\text{$A$}}
\end{picture}
\caption{Four trapezoids\label{ft}}
\end{figure}

\begin{figure}
\unitlength=1.2mm
\begin{picture}(100, 70)(10,10)

\put(20,40){\line(1, 2){10}}
\put(20,40){\line(1,-2){10}}
\put(30,60){\line(1, 0){70}} 
\put(30,20){\line(1, 0){70}} 
\put(90,40){\line(1, 2){10}}
\put(90,40){\line(1,-2){10}}

\put(20,40){\circle*{1}}
\put(90,40){\circle*{1}}
\put(100,60){\circle*{1}}
\put(100,20){\circle*{1}}

\put(15,40){\text{$p$}}
\put(95,40){\text{$r$}}

\put(105,60){\text{$q_2$}}
\put(105,20){\text{$q_1$}}

\put(20,40){\line(4,1){80}}
\put(20,40){\line(4,-1){80}}
\put(40,40){\text{$\th$}}

\end{picture}
\caption{Chevron\label{c}}
\end{figure}

Rotate the chevron around the point $p$ by the angle $\th$
thus taking the point $q_1$ to the point $q_2$. Let $q_3$ be the
image of $q_2$ under this rotation.  Continue rotating the chevron
around $p$ to obtain the sequence $q_1,q_2,q_3,q_4,\ldots$ of
all equally spaced points on the circle centered at $p$ and with
radius $d(p,q_1)$.

Now fix $\epsilon>0$ and consider any continuous transformation of the plane
which keeps the points of the new trapezoid within $\epsilon$ of the old but
changes the angle $\th$.  For example, we could slide the top $B$ a little
to left or right keeping  the bottom $D$ fixed.  This would change the angle
$\th$ continuously say from $\th_0$ to some slightly large value
$\th_1$. To draw the diagram imagine
that $0<\th_1-\th_0<(.0001)\th_0$.  Fix $n$ such 
$n(\th_1-\th_0)>2\pi$.  Note that the sequence 
$$q_1^{\th_1},q_2^{\th_1},q_3^{\th_1},\ldots,q_n^{\th_1}$$
rotates around $p$ at least one more time than the sequence
$$q_1^{\th_0},q_2^{\th_0},q_3^{\th_0},\ldots,q_n^{\th_0}.$$
So by the intermediate value theorem there is a $\th$ with
$\th_0\leq\th\leq\th_1$ with $q_n^\th=q_0^\th$ and another
where $q_n^\th$ is on the opposite side of $p$ to $q_0^\th$.
This implies that there exists $\th$ with $\th_0\leq\th\leq\th_1$
such that the $d(q_1^\th,q_n^\th)$ is exactly the length of say the side 
$B$. Let $T^\prime$ be this trapezoid.

Now we prove that set of vertices of $T^\prime$ 
is Jackson.  Suppose for contraction that $S\sq\reals^2$
meets every isomorphic copy of the vertices of $T^\prime$ in 
exactly one point. 

Suppose $p\in S$.  Then we know that the circle with center
$p$ and radius $d(p,r)$ cannot be a subset $S$. Hence by rotating $S$
we may assume that $r\notin S$.   Since each of
the rightmost two trapezoids must meet $S$, it must be that both
$q_1$ and $q_2$ are in $S$.  If rotate the chevron 
around $p$ by $\theta$ we note that since $q_2\in S$ it must
be that $q_3\in S$.  Similarly it
must be that all $q_i\in S$ in particular $q_n\in S$.  This is
a contradiction since $d(q_1,q_n)$ is the length of $B$ so there
is an isomorphic copy of $T^\prime$ containing two points of $S$.

Other distortions of the trapezoid $T$ which would work as well
to obtain $T^\prime$ would be to 
\begin{itemize}
\item change the height while leaving $B$ and $D$ the same length or 
\item taking any of the vertices and move it slightly to right while
fixing the other vertices.
\end{itemize}
In short anything which continuously changes the angle $\theta$.
For example, if $T$ is a parallelogram, then $T^\prime$ can be
taken to be a parallelogram with the same length sides but
a slight change in the interior angles.
\qed

This proposition suggests the obvious question of whether every four point set
is within $\epsilon$ of a Jackson set.  While the above argument shows that
every four point collinear set is $\epsilon$-close to a trapezoid, we did not
know whether every four point collinear set is $\epsilon$-close to a four
point collinear set which is Jackson.   This question was answered by Louis
Leung:

\begin{prop} (Leung) Every 4-point collinear set is $\epsilon$-close to a
4-point collinear set which is Jackson. \end{prop}

\proof

Let $\{A,B,C,D\}$ be a 4-point set on a straight line \textit{l}.  Reflect the
line across the midpoint of  the line segment $\overline{CD}$. We let $A', B'$
be image of $A,B$ on \textit{l}, hence $B'-D=C-B$ and $A'-B'=B-A$.   Note that
by our usual reflection argument, if $A$ is in a Steinhaus set $S$ then exactly
one of $B'$ or $A'$ lies in $S$.

\def\SS{{\mathcal C}}

We consider two circles $\SS_{1}$ and $\SS_{2}$ centered at $A$ whose radii are 
$d(A,A')$ and $d(A,B')$, respectively.  Let $\alpha, \beta$ be the angles
subtended by a chord of length $d(D,C)$ in $\SS_{1}$ and $\SS_{2}$,
respectively.

Suppose $\epsilon>0$ is fixed.  Without loss of generality we may assume
$\epsilon<d(B,C)$.  Let $C''$ be a point between $B$ and $C$ and $D''$ a point
between $D$ and $B'$ so that $d(C,C'')=d(D,D'')=\delta$
where $ 0<\delta<\epsilon$, and  $\alpha'',\beta''$ be the angles subtended by
a chord of length  $d(D'',C'')$ in $\SS_{1}$ and $\SS_{2}$, respectively. 
Since $d(C,D)<d(C'',D'')$, $\alpha''+\beta''>\alpha+\beta$ and there is
$n\in\mathbb{N}$ such that $n(\alpha''+\beta''-(\alpha+\beta))>2\pi$.  By the
intermediate value theorem, there is a positive $0<\kappa<\epsilon$ such that
if $\tilde{C}$ is between $C$ and $C''$, $\tilde{D}$ is between $D$ and $D''$,
and $d(\tilde{D},D)=d(\tilde{C},C)=\kappa$, then the corresponding angles
$\tilde{\alpha}$ and $\tilde{\beta}$ (the angles subtended by a chord of length
$d(\tilde{C},\tilde{D})$ in $\SS_{1}$ and $\SS_{2}$, respectively) satisfy 

$$n(\tilde{\alpha}+\tilde{\beta})=2m\pi+\theta$$

\noindent for some $m,n\in\mathbb{N}$ where $\theta$ is the angle subtended by
a chord of length $d(A,B)$ in $\SS_{1}$.

\par\bigskip {\bf Claim}
$\{A,B,\tilde{C},\tilde{D}\}$ is Jackson.

\proof Suppose not.  Let $S$ be a Steinhaus set for
$\{A,B,\tilde{C},\tilde{D}\}$.  By translating the set if necessary, we may
assume $A\in S$.  Let $A'$, $B'$, $\SS_{1}$ and $\SS_{2}$ be defined as above. 
Note $\{\tilde{C}, \tilde{D},B',A'\}$ is an isometric copy of
$\{A,B,\tilde{C},\tilde{D}\}$. By our definition of $S$, if $x$ and $y$ lie on
the same half line originating at A such that $x\in \SS_{1}$ and $y\in \SS_{2}$,
either x or y (but not both) must be a member of $S$.  (We may think of, by
rotation if necessary, $x$ as $A'$ and $y$ as $B'$.) Note that $S$ must contain
at least one point in $\SS_{1}$, for otherwise $\SS_{2}\subset S$ and $\SS_{2}$ has a
chord of length $d(A,B)$, a contradiction.

Therefore we may assume $A'\in S$.  We let $A^{'}_{0}=A'$ and $A^{'}_{n+1}$ be
the image of $A^{'}_{n}$ after rotation by $\tilde{\alpha}+\tilde{\beta}$ with
respect to A.  By our definition of $\tilde{\alpha}$, $\tilde{\beta}$ and $S$,
We know $A^{'}_{m}\in S$ for all $m\in\mathbb{N}$.  By letting $m=n$ where $n$
is as given in equation (1), we know that S contains 2 points a distance of
$d(A,B)$ apart, contradicting our choice of $S$. 

Therefore $\{A,B,\tilde{C},\tilde{D}\}$  is Jackson.
\qed

The following still seems to be open:

\begin{question}
Is every four point set in $\reals$
$\epsilon$-close to a four point subset of $\reals$
which is Jackson?
\end{question}

\end{document}